\documentclass[12pt]{article}
\usepackage{amsfonts}
 \usepackage{amsthm}
 \usepackage{amsmath}
 \usepackage{amscd}
 \usepackage{amssymb}
 \usepackage{latexsym}
 \usepackage{epsfig}

\pdfoutput=1
\usepackage[ascii]{inputenc}
\usepackage[T1]{fontenc}
\usepackage{amsmath,amssymb,amsfonts,textcomp,epsfig}
\usepackage{hyperref}
\usepackage{color}
\usepackage{calc}

\graphicspath{{figs/}}

\def\d{{\rm d}}

\def\C{{\mathbb C}}

\def\Z{{\mathbb Z}}
\def\R{{\mathbb R}}

\newtheorem{theorem}{Theorem} 
\newtheorem{lemma}{Lemma}
\newtheorem{proposition}{Proposition}
\newtheorem{corollary}{Corollary}

\title{Perturbations of quadratic Hamiltonian two-saddle cycles.}
 \author{Lubomir Gavrilov \\
 \normalsize \it Institut de Math\'{e}matiques de Toulouse, UMR 5219\\
 \normalsize \it Universit\'{e} de Toulouse, CNRS\\
   \normalsize \it UPS IMT, F-31062 Toulouse Cedex 9, \normalsize \it   France
     \\ Iliya D. Iliev\\
\normalsize \it Institute of Mathematics, Bulgarian Academy of Sciences\\
\normalsize \it Bl. 8, 1113 Sofia, Bulgaria  }
\begin{document}
\maketitle
\begin{abstract}
We prove that the number of limit cycles, which bifurcate
from a two-saddle loop of a planar quadratic Hamiltonian system, under an
arbitrary quadratic  deformation, is less than or equal to three.
\end{abstract}
\newpage
\tableofcontents
\section{Introduction}
The theory of plane polynomial quadratic differential systems
\begin{equation}
 \left\{
\begin{array}{lcr}
   \dot{x}   &= &P(x,y) \\
       \dot{y}   &= &Q(x,y)
\end{array} \right.
\end{equation}
is one of the most classical branches of the theory of two-dimensional autonomous
systems. Despite of the great theoretical interest in studying of such systems,
few is known on their qualitative properties. Let $H(2)$ be the maximal number
of limit cycles, which such a system can have. It is still not known whether $H(2)< \infty$
(or $H(k) < \infty$ for a polynomial system of degree $k$). A survey on the
state of art until 1966 was given by Coppel \cite{copp66} where some basic
and specific properties of the quadratic systems are discussed.

It was believed for a long time that $H(2)=3$, see e.g. \cite{pela58}, until Shi
Song Ling gave in 1980 his famous example of a quadratic system with four limit
cycles \cite{shi80}.

In 1986 Roussarie \cite{rous86b} proposed a local approach to the global
conjecture $H(k)<\infty$, based on the observation that if the
cyclicity is infinite, then a limiting periodic set will exist with infinite
cyclicity. All possible 121 limiting periodic sets of quadratic systems were
later classified in \cite{drr94}.

Of course, it is of interest to compute explicitly the cyclicity of concrete
limiting period sets, the simplest one being the equilibrium  point. It is another
classical result, due to Bautin (1939), which claims that the cyclicity of a
singular point of a quadratic system is at most three. The cyclicity of
Hamiltonian quadratic homoclinic loops is two  \cite{hoil94a,ilie96}, and for
the reversible ones see \cite{heli95}.

In his controversial paper \cite{zola95}, \.{Z}o\l\c{a}dek proved that the cyclicity
of the Melnikov functions near quadratic triangles (three-saddle loops) or
segments (two-saddle loops) is respectively three and two. From this he deduced
that the cyclicity of the triangle or the segment itself is also equal to three
or two, respectively. As we know now, this conclusion is not always true.
Namely, in the perturbed Hamiltonian case, not all limit cycles near a
polycycle are "shadowed" by a zero of a Melnikov function. The bifurcation of
"alien" limit cycles is a new phenomenon discovered recently by Caubergh, Dumortier
and Roussarie \cite{cdr05, duro06}. Li and Roussarie \cite{liro04} later
computed the cyclicity of quadratic Hamiltonian two-loops, when they are
perturbed "in a Hamiltonian direction". In the case of a more general
perturbation they only noted that "some new approach may be needed". 

One of the most interesting developments in this field, starting from the
series of papers by Petrovskii and Landis \cite{pela58}, is the proliferation
of complex methods, as it can be seen from the 2002 survey of Ilyashenko
\cite{ilya02}. A particular interest is given to the study of different
infinitesimal versions of the 16th Hilbert problem.
Thus, G.S. Petrov \cite{petr88a} used the argument principle to evaluate the
zeros of suitable complete Abelian integrals, which on its turn produces an
upper bound for the number of limit cycles, which a perturbed quadratic system
of the form
$$
 \left\{
\begin{array}{lcr}
   \dot{x}   &= &y + \varepsilon P(x,y) \\
       \dot{y}   &= &x-x^2 + \varepsilon Q(x,y)
\end{array} \right.
$$
may have. The result was later generalized for the perturbations of arbitrary
generic cubic Hamiltonians in \cite{hoil94b,gavr01}.

The present paper studies the cyclicity of quadratic Hamiltonian monodromic
two-loops, as on Fig. \ref{fig1}. We use complex methods, in the spirit of
\cite{gavr11, gavr11a}, which can also be seen as a far going generalization of
the original Petrov method. Our main result is  that at most three limit cycles
can bifurcate from such a two-loop (Theorem \ref{main1}), although we did not
succeed to prove that this bound is exact. It is interesting to note, that even for a generic
quadratic perturbation, two limit cycles can appear near a two-saddle loop, while at the same
time the (first) Poincar\'e-Pontryagin (or Melnikov) function exhibits only one zero.
The appearance of the missing alien limit cycle is discussed in the Appendix.

Our semi-local results, combined with the known cyclicity of open period annuli lead also to some global results,
formulated in Section \ref{globalsection}.

\section{Statement of the result}
\begin{figure}
\begin{center}
 \def\svgwidth{16cm}
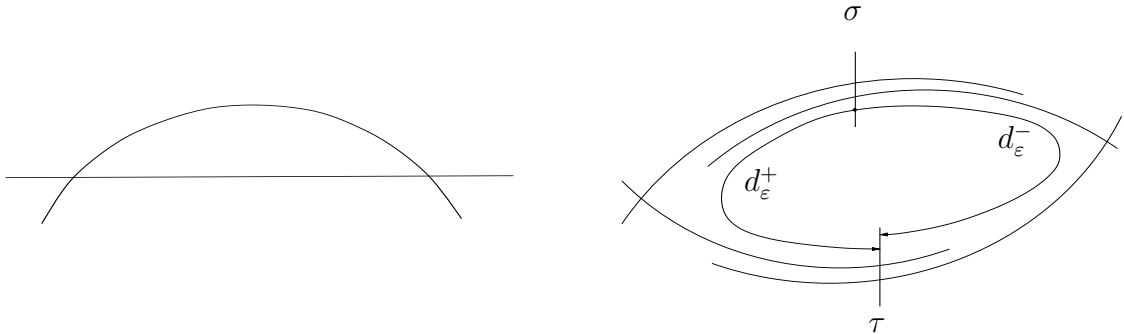
\end{center}
\caption{Monodromic two-saddle loop and the Dulac maps $\d_\varepsilon^\pm$}
\label{fig1}
\end{figure}

Let $X_\lambda$, $\lambda\in \R^{12}$ be the (vector) space of all quadratic
planar vector fields, and let $X_{\lambda_0}$ be a planar quadratic vector field
which has two non-degenerate saddle points $S_1(\lambda_0), S_2(\lambda_0)$
connected by two
heteroclinic connections $\Gamma_1, \Gamma_2$, which form a monodromic two-loop
as on Figure \ref{fig1}. The  union $\Gamma= \Gamma_1\cup \Gamma_2$ will be
referred to as a non-degenerate two-saddle loop.
The cyclicity $Cycl(\Gamma, X_\lambda)$ of the two-saddle loop $\Gamma$ with
respect to the deformation $X_\lambda$ is the maximal number of limit cycles
which $X_\lambda$ can have in an arbitrarily small neighborhood of $\Gamma$,
 as $\lambda$ tends to $\lambda_0$, see \cite{rous98}.

 In the present paper we shall suppose in addition, that $X_{\lambda_0}$ is a
 Hamiltonian vector field
 \begin{equation}
\label{hamiltonian}
X_{\lambda_0}=X_H: \left\{
\begin{array}{lcr}
   \dot{x}   &= &H_y \\
       \dot{y}   &= & -H_x
\end{array} \right.
\end{equation}
where $H$ is a bivariate polynomial of degree three.
Our main result is the following

\begin{theorem}
\label{main1}
The cyclicity of every non-degenerate Hamiltonian two-saddle loop, under an
arbitrary quadratic deformation, is at most equal to three.
\end{theorem}
The result will be proved by making use of complex methods, as explained in
\cite{gavr11}, combined with the precise computation of the so called higher
order Poincar\'e-Pontryagin (or Melnikov) functions, which can be found in
\cite{ilie98}.

\subsection{Outline of the proof of Theorem 1}
\subsubsection{Principalization of the Bautin ideal}
Let $h \mapsto P_\lambda(h)$ be the first return map associated to the deformed
vector field $X_\lambda$ and the period annulus of $X_{\lambda_0}$, bounded by
$\Gamma$. Consider the Bautin ideal
$${\cal B}= \langle a_k(\lambda)\rangle \subset \C [\lambda]$$
 generated by the coefficients of the expansion
$$
P_\lambda(h)-h = \sum_{k=1}^\infty a_k(\lambda) h^k  .
$$
In the quadratic case under consideration its computation is well known, and
goes back to Bautin, see \cite{ilya08} for details. It  has three generators,
in particular the ideal  is not principal. By making use of the Hironaka's
desingularization theorem, we can always assume that ${\cal B}$ is
"locally principal". Namely, by abuse of notation, let ${\cal B}$ be the ideal
sheaf generated by the Bautin ideal, in the sheaf of analytic functions
${\cal O}_X$ on $X$. The parameter space $X=\R^{12}$ can be replaced by a new
smooth real analytic variety $\tilde{X} $, together with a proper analytic map
$$
\pi: \tilde{X} \rightarrow X
$$
such that the pull back  $\pi^*{\cal B}$ is a  principal ideal sheaf. This
means that for every point $\tilde\lambda \in \tilde{X}$ there is a neighborhood
$U$, such that the ideal $\pi^*{\cal B}(U)$ of the ring
${\cal O}_{ \tilde{X} }(U)$ is a principal ideal,
see \cite[section 2.1]{gavr08} and Roussarie \cite{rous01}.

The cyclicity at a point $\lambda_0 \in X$ is the lower upper bound of the
cyclicities computed at points of the compact set $\pi^{-1}(\lambda_0)$. As
the cyclicity is an upper semi-continuous function in $\tilde{\lambda}_0 \in
\pi^{-1}(\lambda_0)$, and $ \pi^{-1}(\lambda_0)$ is compact, then there is a
$\tilde{\lambda}_0 \in \pi^{-1}(\lambda_0)$ at which the cyclicity
$Cycl(\Gamma, X_{\tilde\lambda})$ is maximal. It suffices therefore to compute
this cyclicity.

In more down to earth terms, the above considerations show that, after
appropriate analytic change of the parameters $\lambda=\lambda(\tilde{\lambda})$,
we can always suppose that the localization of the Bautin ideal at $\lambda_0$ is
a principal ideal of the ring of germs of analytic functions at $\lambda_0$.
We denote its generator (according to the tradition) by $\varepsilon$.
The power series expansion of the first return map takes therefore the form
\begin{equation}
\label{melnikov}
P_\lambda(h)= h + \varepsilon^k [M_k(h) + O(\varepsilon)],
\quad M_k\neq 0
\end{equation}
where $M_k$ is the $k$-th order Melnikov function, associated to $P_\lambda$.
The function $O(\varepsilon)$, by abuse of
notation, depends on $h, \lambda$ too, but it is of $O(\varepsilon)$ type
uniformly in $h, \lambda $, where $h$ belongs to a compact complex domain in
which the return map is regular. The principality of the Bautin ideal is
equivalent to the claim, that $M_k(h)$ is not identically zero.
The perturbed Hamiltonian vector field $X_{\varepsilon,\lambda}$ can be
supposed on its turn of the form
 \begin{equation}
\label{perturbedfield}
X_{\varepsilon,\lambda}: \left\{
\begin{array}{lcr}
   \dot{x}   &= &H_y + \varepsilon Q(x,y,\lambda,\varepsilon) \\
       \dot{y}   &= & -H_x - \varepsilon P(x,y,\lambda,\varepsilon)
\end{array} \right.
\end{equation}
where $P,Q$ are quadratic polynomials in $x,y$ with coefficients depending
analytically in $\varepsilon, \lambda$.
Of course, we shall need an explicit expression for $M_k(h)$ which depends
also on the unknown parameter value $\tilde{\lambda}_0 \in \pi^{-1}(\lambda_0)$.
Taking analytic curves
\begin{equation} \label{analytic}
\varepsilon \mapsto \lambda(\varepsilon), \quad \lambda(0)=\lambda_0
\end{equation}
 we get from (\ref{melnikov})   
$$
P_{\lambda(\varepsilon)}(h)= h + \varepsilon^k [M_k(h) + O(\varepsilon)],
\quad M_k\neq 0
$$
which allows one to compute $M_k$ by only making use of analytic one-parameter
deformations and the Fran\c{c}oise algorithm \cite{fran96}.
The general form of the first non-vanishing Melnikov function with respect to
any analytic curve of the form (\ref{analytic})
 in the Hamiltonian (or more generally, integrable) quadratic case is computed
 in \cite{ilie98}.

By abuse of notation, from now on, the return map of the form (\ref{melnikov}),
will be denoted by $P_\varepsilon$, where $\varepsilon$ is the generator of
the localized Bautin ideal.

\subsubsection{The Petrov trick and the Dulac map}
The limit cycles of $X_\lambda$ are the fixed points of $P_\varepsilon$.
We are going to study these fixed points in a complex domain, where they
correspond to complex limit cycles. $P_\varepsilon$ is obviously a composition
of two Dulac maps $d^\pm(\varepsilon)$ as on Fig. \ref{fig1}
$$
P_\varepsilon= (d_\varepsilon^-)^{-1}\circ d_\varepsilon^+
$$
so the fixed points $h$ of $P_\varepsilon$ are the zeros of the displacement
map $d_\varepsilon^+ - d_\varepsilon^-$. In a complex domain this map has two
singular points corresponding to the saddles $S_\pm(\varepsilon)$ and we shall
study its zeros in the complex domain
${\cal D}_\varepsilon$, shown on Fig. \ref{fig5}. This domain is bounded by a
circle, by the segment $(S_+(\varepsilon), S_-(\varepsilon))$, and by the zero
locus of the imaginary part of $d^+_\varepsilon$. The number of the zeros of
$d_\varepsilon^+ - d_\varepsilon^-$ in ${\cal D}_\varepsilon$ is computed
according to the argument principle: it equals the increase of the argument
along the boundary of ${\cal D}_\varepsilon$.

Along the circle and far from the critical points, the displacement function is
"well" approximated by $\varepsilon^k M_k(h)$ which allows one to estimate the
increase of the argument.

Along the segment $(S_+(\varepsilon), S_-(\varepsilon))$ the zeros of the
imaginary part of the displacement function coincide with the fixed points of
the holomorphic holonomy map along the separatrix through $S_-(\varepsilon)$.
The zeros are therefore well approximated, similarly to (\ref{melnikov}),
by an Abelian integral along the cycle $\delta_-(h)$ in the fibers of $H$,
vanishing at $S_-(0)$. This observation may be seen as a far going
generalization of the so called Petrov trick, see \cite{gavr11a} for details.

Along the zero locus of the imaginary part of $d^+_\varepsilon$, the zeros of
imaginary part of the displacement map coincide with the fixed points of the
composition of the holonomies associated to the separatrices through
$S_-(\varepsilon)$ and $S_+(\varepsilon)$. As this map is holomorphic,
it is similarly approximated by the zeros of an Abelian integral along
$\delta_-(h)+ \delta_+(h)$, where $\delta_\pm(h)$ are cycles in the fibers of
$H$, vanishing at $S_\pm(0)$ respectively.

Thus, to count the number of the limit cycles, it is enough to inspect the
behavior of certain Abelian integrals.

\section{Abelian integrals related to quadratic perturbations of reversible
quadratic Hamiltonian vector fields.}
\label{abelianintegrals}
In this section we recall the Abelian integrals, involved in the proof of
Theorem \ref{main1}, and establish their properties. The details can be found
in \cite{ilie96,ilie98}.

Consider the  quadratic reversible Hamiltonian system $dH=0,$
where the Hamiltonian function is taken in the normal form \cite{ilie96}
\begin{equation}
\label{ham}
H(x,y)=x[y^2+ax^2-3(a-1)x+3(a-2)], \quad a  \in  {\mathbb R}\end{equation}
The Hamiltonian system has a center  $C_0=(1,0)$ on the level set $H=t_0=a-3$.
It is surrounded by a saddle connection containing two saddles
$S_ \pm= (0,  \pm\sqrt{3(2-a)})$
if and only if the parameter $a$ takes values in $(-1,2)$. This connection
is a part of the zero-level set $H=t_s=0$.
When $a\in(0,2)$, there is a second center
$$C_1=\left(\frac{a-2}{a} ,0\right), \quad H(C_1)=t_1=
\frac{(a+1)(a-2)^2}{a^2}$$
surrounded by other part of the zero level set and containing
the same two saddles.

Let $\delta(t) \subset \{H=t\}$ be a continuous family of ovals
surrounding a center.
Take a small quadratic one-parameter perturbation
\begin{equation}
\label{heps}
dH+\varepsilon\omega=0, \quad \omega=\omega(\varepsilon)=f(x,y,\varepsilon)dx
+g(x,y,\varepsilon)dy,
\end{equation}
where $f$, $g$ are real quadratic polynomials
of $x,y$ with coefficients analytic with respect to the small parameter
$\varepsilon$.
Then the first return map  $P_\varepsilon$ near an oval $\delta(t)$ is
well defined and has the form
\begin{equation}
\label{returnmap}
P_\varepsilon(t)=t + \varepsilon M_1(t)+\varepsilon^2M_2(t)
+\varepsilon^3M_3(t)+\ldots,
\end{equation}
One may show, by making use of \cite[Theorem 2]{gail05}, that the first
non-vanishing Poincar\'e-Pontryagin-Melnikov function $M_k$ associated to an
arbitrary polynomial perturbation is an Abelian integral. More precisely,
we have
\begin{theorem}
[\cite{ilie98}]
In the quadratic case $M_k$ takes the form
\begin{equation}
M_1(t)=\int_{\delta(t)}[  \alpha+\beta x]ydx,\quad M_k(t)=
\int_{\delta(t)}[  \alpha+\beta x+\gamma x^{-1}] ydx,\quad k\geq 2,
\label{mk}
\end{equation}
where $ \alpha,\beta,\gamma $ are appropriate constants depending on
the perturbation.
\end{theorem}
Consider the Abelian integrals
$$J_k(t)=\int_{\delta(t)} x^kydx,  \quad k\in {\mathbb Z}$$
(oriented clockwise - along with the Hamiltonian vector field).
\begin{lemma}[\cite{ilie96}]
\label{pfs}
The integrals $J_k(t)$,   $k=-1,0,1$  satisfy the following system with respect
to $t$:
\begin{equation}
\label{p-f}
\begin{array}{l}
tJ_{-1}'+(4-2a)J_0'+(a-1)J_1'=\frac 13 J_{-1},\\
(1-a)tJ_{-1}'+2atJ_0'+(3+2a-a^2)J_1'=\frac 43a J_0,\\
(a-2)tJ_{-1}'+(2-2a)tJ_0'+at J_1'= \frac 32(1-a) J_0 + aJ_1.
\end{array}
\end{equation}
\end{lemma}

\begin{lemma}
\label{asymp}
{\it The integrals $J_k(t)$,   $k=-1,0,1$  have the following asymptotic
expansions near $t=-0$:}
\begin{equation}
\label{exp}
\begin{array}{rl}
J_{-1}(t)&=-2\sqrt{3(2-a)}[1-\frac{a-1}{12(a-2)^2}t-\frac{11a^2-22a+15}
{576(a-2)^4}t^2-\frac{35(a-1)(5a^2-10a+9)}{20736(a-2)^6}t^3+\ldots ]\ln t\\[2mm]
&+ a_0+a_1t+a_2t^2\ldots,\\[2mm]
J_0(t)&=-2\sqrt{3(2-a)}[-\frac{1}{6(a-2)}t-\frac{a-1}{48(a-2)^3}t^2
-\frac{85a^2-170a+105}{10368(a-2)^5}t^3+\ldots ]\ln t\\[2mm]
&+ b_0+b_1t+b_2t^2+\ldots,\\[2mm]
J_1(t)&=-2\sqrt{3(2-a)}[-\frac{1}{72(a-2)^2}t^2
-\frac{5(a-1)}{864(a-2)^4}t^3+\ldots]\ln t+ c_0+c_1t+c_2t^2\ldots.
\end{array}
\end{equation}
\end{lemma}
This lemma is a consequence of the following basic property of system
(\ref{p-f}):

\begin{lemma}
\label{asymp2}
{\it If $a\neq 0$, a fundamental system of solutions $(J_{-1}, J_0, J_1)^t $ of
$(\ref{p-f})$ near $t=0$ is the following:}
\begin{equation}
\label{fund}
\begin{array}{l}
P(t)=\left(\begin{array}{c}
3(a-1)\\   \frac{3(3+2a-a^2)}{4a}\\
\frac{9(a-1)(3+2a-a^2)}{8a^2}\end{array}\right)+\left(\begin{array}{c}
0\\0\\1\end{array}\right)t,\\[3mm]
Q(t)=\left(\begin{array}{c}1\\0\\0\end{array}\right)-
\left(\begin{array}{c}\frac{a-1}{12(a-2)^2}\\\frac{1}{6(a-2)}\\0
\end{array}\right)t-
\left(\begin{array}{c}\frac{11a^2-22a+15}{576(a-2)^4}\\
\frac{a-1}{48(a-2)^3}\\ \frac{1}{72(a-2)^2}
\end{array}\right)t^2-
\left(\begin{array}{c}\frac{35(a-1)(5a^2-10a+9)}{20736(a-2)^6}\\
\frac{85a^2-170a+105}{10368(a-2)^5}\\ \frac{5(a-1)}{864(a-2)^4}
\end{array}\right)t^3+\ldots,\\[3mm]
R(t)=Q(t)\ln t+S(t),\end{array}
\end{equation}
{\it with $S(t)$ analytic function in a neighbourhood of $t=0$}.
\end{lemma}

\vspace{2ex}
\noindent
{\bf Proof.} Rewrite system (\ref{p-f}) in the form $(A_1t+A_0)J'=BJ$.
System (\ref{p-f}) has at its critical value $t=0$ a triple
characteristic exponent equal to zero, while its characteristic exponents at
infinity are $-\frac13$, $-\frac23$, $-1$. Hence, there is a polynomial
solution $P(t)$ of degree one which is easy to find. To calculate $Q(t)$,
we replace
$$Q(t)=q_0+q_1t+q_2t^2+q_3t^3+q_4t^4+\ldots$$
in the system to obtain recursive equations
\begin{equation}
\label{rec}
(jA_1-B)q_j+(j+1)A_0q_{j+1}=0,\quad j=0, 1, 2, 3, \ldots.
\end{equation}
The third equation in the system obtained for $j=0$ implies that
$\frac32(1-a)q_{0,0}+aq_{0,1}=0$ where $q_0=(q_{0,-1}, q_{0,0}, q_{0,1})^\top$.
Therefore, one can choose without loss of generality $q_0=(1,0,0)^\top$.
Then any analytic solution of (\ref{p-f}) would be a unique linear combination
of $P(t)$ and $Q(t)$. Fixing in such a way $Q(0)$, then $q_1$, $q_2$ and so on
are uniquely determined from system (\ref{rec}).

Now, if we take a linear combination $\tilde Q$ of $P$ and $Q$ and replace
$R(t)=\tilde Q(t)\ln t+S(t)$ in the system $(A_1t+A_0)R'(t)=BR(t)$,
we obtain $A_1\tilde Q+t^{-1}A_0\tilde Q+(A_1t+A_0)S'=BS$. Hence,
$A_0\tilde Q(0)=0$ which means that $\tilde Q(t)$ is proportional to $Q(t)$.
Therefore one can simply take $\tilde Q=Q$. $\Box$

\vspace{2ex}
\noindent
{\bf Proof of Lemma \ref{asymp}.} Let $x_1$ be the (smaller) positive root of
the equation $r(x)=-ax^2+3(a-1)x-3(a-2)=0$ where $a\in(-1,2)$.
Then, $J_k(0)=\int_{\delta(0)}x^kydx=2\int_0^{x_1}x^k\sqrt{r(x)}dx$ for
$k=0,1$. Therefore
$$\frac32(a-1)J_0(0)-aJ_1(0)=\int_0^{x_1}r'(x)\sqrt{r(x)}dx=
-\frac23(3(2-a))^{3/2}.$$
On the other hand, the third equation of (\ref{p-f}) implies
$$ \frac32(1-a)J_0(0)+aJ_1(0)=(a-2)\left.(tJ_{-1}'(t))\right|_{t=0}.$$
Finally, if $J(t)=\lambda[Q(t)\ln t +S(t)]+\mu P(t)+\nu Q(t)$,
then $\left. (tJ_{-1}'(t))\right|_{t=0} =\lambda=-2\sqrt{3(2-a)}$.
Case $a=0$ follows by continuity. $\Box$

%

\section{Cyclicity of  two-saddle cycles}
In the section we prove Theorem \ref{main1}.

We shall prove it in several steps. A plane quadratic Hamiltonian system with a
two-saddle loop can be written, up to an affine change of the variables, in the
form $dH=0$ where $H$  is in the form (\ref{ham}).
\subsection{The case $M_1\neq0$}
\label{m10}
In this section we consider the  perturbed quadratic plane quadratic
Hamiltonian system (\ref{heps}) under the generic assumption that
$$
M_1(t)= \int_{\delta(t)} \omega|_{\varepsilon=0} =
\int\int_{\{H\leq t\}}[  \alpha+\beta x]dx dy
$$
is not identically zero.

Due to Lemma \ref{asymp}, $M_1(t)$ vanishes identically in a co-dimension two
analytic set defined by $\{\alpha=\beta=0\}$.
The Poincar\'e-Pontryagin function $M_1$ is well defined at $t=0$ in which case
it is the well known Melnikov integral along the heteroclinic loop $\delta(0)$.
It is classically known that when $M_1(t)\not\equiv 0$, the vanishing of the
Melnikov integral $M_1(0)$ is a necessary condition for a bifurcation of a
limit cycle (and in the opposite case the heteroclinic loop is broken under the
perturbation)
\begin{proposition}
\label{nocycles}
If $M_1(0) \neq 0$, then no limit cycles bifurcate from the two-saddle loop
$\Gamma$.
\end{proposition}
{\bf Proof.} Suppose that there is a sequence of limit cycles
$\delta_{\varepsilon_i}$ of (\ref{heps}) which tend to $\Gamma$ as
$\varepsilon_i$ tends to $0$. Then
$$
0= -\int_{\delta_{\varepsilon_i}} dH= \varepsilon_i
\int_{\delta_{\varepsilon_i}} \omega
$$
which implies
$$
  0= \lim_{\varepsilon_i \rightarrow 0}  \int_{\delta_{\varepsilon_i}} \omega =
  \int_{\Gamma} \omega|_{\varepsilon=0} = M_1(0)  .
$$
$\Box$

The complete Abelian integral  $M_1(t)$  has the following convergent expansion
near the critical saddle value $t=0$
\begin{equation}
\label{expansion}
M_1(t)= d_0+ d_1 t \ln t + d_2 t + d_3 t^2 \ln t + \dots
\end{equation}
Let $\delta_\pm(t)\in H_1(\Gamma_t,\Z)$,
$\Gamma_t= \{(x,y)\in \C^2: H(x,y)=t\}$,
be the two continuous families of cycles, vanishing at the saddle points
$S_\pm$ respectively, with  orientations chosen in a way that for the
respective intersection indices holds
\begin{equation}
\label{orientation}
\delta\cdot \delta_+ = \delta\cdot \delta_- = -1 .
\end{equation}
Then
\begin{equation}
\label{m1}
M_1(t)=\int_{\delta(t)} \omega_0 = \frac{ \ln t}{2\pi \sqrt{-1} }
(\int_{\delta_+(t)} \omega_0 + \int_{\delta_-(t)} \omega_0)+ d_0+ d_2 t+ O(t^2)
\end{equation}
where $\omega_0= \omega|_{\varepsilon=0}$. The involution
$(x,y)\rightarrow (x,-y)$ leaves the level set $\{H=h\}$ invariant, reversing
its orientation. Therefore it acts on $\delta, \delta_\pm$ as follows
$$
\delta \rightarrow - \delta, \;\; \delta_- \rightarrow -\delta_+, \;\;
\delta_+ \rightarrow -\delta_- ,
$$
which implies
\begin{equation}
\label{m1m1}
\int_{\delta_+(t)} \omega_0 = \int_{\delta_-(t)} \omega_0 .
\end{equation}
Let $h^\varepsilon_{\delta_\pm}$ be the two holonomy maps associated to the
separatrices of the perturbed foliation, intersecting the cross-section
$\sigma$. There are two-possible orientations for the loop defining the
holonomy, this  corresponds  to a choice of orientation of $\delta_\pm$,
see (\ref{orientation}). Similarly to (\ref{returnmap}) we have
\begin{eqnarray}
h^\varepsilon_{\delta_+}(t) & = &t + \varepsilon \int_{\delta_+(t)}
\omega_0 + O(\varepsilon^2) \label{hez+}\\
h^\varepsilon_{\delta_-}(t) & = &t + \varepsilon \int_{\delta_-(t)}
\omega_0 + O(\varepsilon^2) \label{hez-}\\
h^\varepsilon_{\delta_+} \circ h^\varepsilon_{\delta_-}(t) & = &
t+\varepsilon (\int_{\delta_+(t)} \omega_0 + \int_{\delta_-(t)}
\omega_0) + O(\varepsilon^2)\\
h^\varepsilon_{\delta_-} \circ h^\varepsilon_{\delta_+}(t) & = &
t+\varepsilon (\int_{\delta_+(t)} \omega_0 + \int_{\delta_-(t)}
\omega_0) + O(\varepsilon^2) \label{hez4}
\end{eqnarray}

\begin{proposition}
\label{d0d1}
If $d_0=d_1=0$, then $\alpha= \beta = 0$.
\end{proposition}
\textbf{Proof.} According to Lemma \ref{asymp} $d_1= \alpha/\sqrt{3(a-2)}$.
If $\alpha=0$ then
$$d_0= \beta J_1(0) \mbox{  where  } J_1(0)= \int\!\!\int_{\{H<0\}}x dx
\wedge dy \neq 0 .\Box
$$

Therefore $M_1\neq 0$ if and only if $|d_0|^2+|d_1|^2\neq 0$, and hence at most
one zero of $M_1$ can bifurcate from $t=0$. Of course, no conclusion about the
number of the limit cycles can be deduced at this stage. For a further use, let
us note that the above implies (see also Proposition \ref{nocycles})
\begin{corollary}
\label{simplezero}
If a limit cycle bifurcates from the two-saddle loop, then
the Abelian integral $\int_{\delta_\pm(t)} \omega_0$ has a simple zero at the
origin.
\end{corollary}

\begin{proposition}
\label{pr5}
If the Melnikov function $M_1$ is not identically zero, then at most two limit
cycles bifurcate from $\Gamma$.
\end{proposition}
\begin{proposition}
\label{pr5+}
There exists a perturbed quadratic system of the form (\ref{perturbedfield})
and $M_1\neq0$, with exactly two limit cycles bifurcating from the two-saddle
loop.
\end{proposition}
The proof of this proposition will be postponed to the Appendix. To the end of
this subsection we shall prove Proposition \ref{pr5}.
 Although our proof will be self-contained, we shall omit some technical
 details, for which we refer to \cite[section 4]{gavr11}.

\begin{figure}
\begin{center}
\resizebox{7cm}{!}
{ 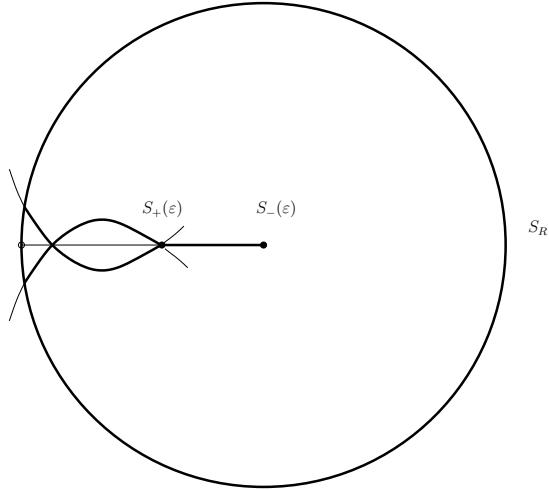}
\end{center} \caption{The domain $\mathcal{D}_\varepsilon$}\label{fig5}
\end{figure}

Consider the Dulac maps $d^+_\varepsilon$,
$d^-_\varepsilon$ associated to the perturbed foliation, and to the cross
sections $\sigma$ and $\tau$, see Fig. \ref{fig1}.
We parameterize each cross-section by the restriction  of the first integral
$f$ on it, and denote $t=f|_\sigma$. Each function $d^\pm_\varepsilon$ is
multivalued and has a critical point at
$S_\pm(\varepsilon)\in \mathbb{R}$, $S_\pm(0)=0$. The  saddle points
$S_+, S_-$ depend analytically on $\varepsilon$. Without loss of generality we
shall suppose that $\varepsilon > 0$ and $S_-(\varepsilon) > S_+(\varepsilon)$,
see Fig. \ref{fig5}.
A limit cycle intersects the cross-section $\sigma$ at $t$ if and only if
$d^+_\varepsilon(t) = d^-_\varepsilon(t)$. Therefore
zeros of the displacement map
$$
d^+_\varepsilon - d^-_\varepsilon = (d^+_\varepsilon \circ
(d^-_\varepsilon)^{-1} - id)\circ d^-_\varepsilon =
(P_\varepsilon - id)\circ d^-_\varepsilon
$$
correspond to limit cycles. Our aim is to bound the number of those zeros
which are real, bigger than $S_-(\varepsilon)$, and tend to $0$ as
$\varepsilon$ tends to $0$. For this, we consider an appropriate complex domain
$\mathcal{D}_\varepsilon$ of the universal covering of
$\mathbb{C}\setminus \{S_+(\varepsilon)\}$ and compute the number of the zeros
of the displacement map, by making use of the argument principle. The reader
may find useful to compare our method, to  the Petrov's method \cite{petr90},
used to compute zeros of complete elliptic integrals. The crucial fact is
that, roughly speaking,  \emph{the monodromy of the Dulac map is the holonomy
of its separatrix}. The analytical counter-part of this statement is that
the zero locus $\mathcal{H}^\pm_\varepsilon$ of the imaginary part of the
Dulac map $d^\pm_\varepsilon$ for $\Re (t) <
S_\pm(\varepsilon)$ is a real-analytic curve in
$\{ \mathbb{R}^2=\mathbb{C} \} \cap \mathcal{D}_\varepsilon$,
defined in terms of the holonomies of the separatrices.
 It follows from  \cite[section 4]{gavr11} that
$$
\mathcal{H}_\varepsilon^+ = \{z\in\mathbb{C}^2 : h^\varepsilon_{\delta_+}(z)
= \overline{z} \},
\mathcal{H}_\varepsilon^- = \{z\in\mathbb{C}^2 :
h^\varepsilon_{\delta_-}(\overline{z}) = z\} .
$$
Note that the above describes, strictly speaking, only one connected component
of $\mathcal{H}_\varepsilon^\pm$, the second one is "complex conjugate" and
defined by a similar formula
$$
\mathcal{H}_\varepsilon^+ = \{z\in\mathbb{C}^2 :
h^\varepsilon_{\delta_+}(\overline{z}) = z \},
\mathcal{H}_\varepsilon^- = \{z\in\mathbb{C}^2 :
h^\varepsilon_{\delta_-}(z) = \overline{z}\} .
$$
By abuse of notation we use $\mathcal{H}_\varepsilon^\pm$ to denote only the
first connected component (the second corresponds to the opposite orientation
of $\delta_\pm$).

The analyticity of the above curves is crucial in computing the complex zeros
of the transcendental Dulac maps. For instance, to compute the number of
intersection points of $\mathcal{H}_\varepsilon^\pm$ with the real axis
$\{z= \bar{z}\}$ we have to solve the equation
\begin{equation}
\label{hez}
h^\varepsilon_{\delta_\pm}(z)=z,
\end{equation}
 and to compute the number of the intersection point of
 $\mathcal{H}_\varepsilon^-$ with $\mathcal{H}_\varepsilon^+$, we have to solve
 the equation
\begin{equation}
\label{hez2}
h^\varepsilon_{\delta_-} \circ h^\varepsilon_{\delta_+}(z) = z .
\end{equation}

Let us define first the complex domain
$\mathcal{D}_\varepsilon$ in which the computation will take place:
it is bounded by the circle
$$
S_R= \{t: |t| = R \},
$$
by the interval $[S_+(\varepsilon),S_-(\varepsilon)]$, and by the zero locus
$\mathcal{H}_\varepsilon^+ $, as it is shown on Fig. \ref{fig5}.

Let $R, \varepsilon_0$ be real numbers subject to certain technical conditions
of the form
$$
1 >> R >> \varepsilon_0 > 0 .
$$
The subsequent computations will hold for all $\varepsilon$, such that
$$
 \varepsilon_0 > \varepsilon  >0 .
$$
We wish to bound the number of the zeros of the displacement map in the domain
$\mathcal{D}_\varepsilon$. If the map were an analytic function in a
neighborhood of the closure of the domain, and non-vanishing on its border,
we could apply the argument principle:
\begin{quote}

\emph{The number of the zeros (counted with multiplicity)  in the complex
domain $\mathcal{D}_\varepsilon$ equals the increment of the argument of this
function along the border of $\mathcal{D}_\varepsilon$, divided by $2\pi$. }
\end{quote}
The above principle holds true with the analyticity condition relaxed: it is
enough that the map allows a continuation on the closure of the domain
$\mathcal{D}_\varepsilon$, considered as a subset of the universal covering
of
$$\mathbb{C}\setminus \{S_+(\varepsilon),S_-(\varepsilon)\} .$$
This is indeed the case, and it remains to assure finally the non-vanishing
property. Along $S_R$ the displacement map  has a known asymptotic behavior and
hence does not vanish. Along the remaining part of the border, including
$S_\pm(\varepsilon)$ the displacement map can have isolated zeros.
For this  we may add to the displacement map a small real constant
$c>0$, sufficiently smaller with respect to $\varepsilon$. The new function
$d^+_\varepsilon - d^-_\varepsilon + c$ which we obtain in this way has at
least so many zeros in $\mathcal{D}_\varepsilon$, as the original
displacement map, but is  non-vanishing on the border of the domain. The
increase of the argument of $d^+_\varepsilon - d^-_\varepsilon + c$ along
$S_R$ will be close to the increase of the argument of
$d^+_\varepsilon - d^-_\varepsilon$ (because $c<<\varepsilon$). At last,
the imaginary parts of $d^+_\varepsilon - d^-_\varepsilon$  and
$d^+_\varepsilon - d^-_\varepsilon+c$ are the same. The intuitive content of
this is that when the displacement map has zeros on the border of the domain,
it will have less zeros in the interior of the domain.

To resume, according to the argument principle, to evaluate the number of the
zeros of the displacement map in the the domain $\mathcal{D}_\varepsilon$, it
is enough to evaluate
\begin{enumerate}
\item The increase of the argument of  the displacement map,  along the circle
$S_R$.
\item The number of the zeros of  the imaginary part of the displacement
 map,  along the  interval $[S_+(\varepsilon),S_-(\varepsilon)]$.
  \item The number of the zeros of  the imaginary part of the displacement map,
  along the real analytic curve $\mathcal{H}_\varepsilon^+$.
 \end{enumerate}
 To the end of the section we evaluate the above quantities.
\begin{enumerate}
  \item By Proposition \ref{nocycles}, if limit cycles bifurcate from the
  double loop, then
$$
d_0= \alpha J_0(0) + \beta J_1(0) = \int\!\!\int_{\{H<0\}}(\alpha + \beta x) dx
\wedge dy  = 0
$$
and hence $\alpha\neq0, \beta \neq 0$. From this we conclude  that the
displacement map along the circle $S_R$ is approximated by  $\varepsilon M_1$
which has as a leading term $t \ln t$ (because $d_0=0$ but $d_1\neq 0$).
The increase of the argument of $t \ln t$, and hence of the displacement map,
along the circle $S_R$ is \emph{close to $2\pi$ but strictly less than $2\pi$}.
  \item The imaginary part of the displacement map,  along the  interval
  $[S_+(\varepsilon),S_-(\varepsilon)]$ equals  the imaginary part of
$d^-_\varepsilon(t)$. Its zeros equal the number of intersection points of
$\mathcal{H}_\varepsilon^+$ with the real axes, which amounts to solve
$h^\varepsilon_{\delta_-}(z)=z$, see
(\ref{hez}). By (\ref{hez-}) the number of the zeros is bounded by the
multiplicity of the holomorphic Abelian integral $\int_{\delta_-(t)}
\omega_0$ having a simple zero at the origin (Corollary \ref{simplezero}).
Note, however, that the holonomy map $h^\varepsilon_{\delta_-}$
has $S_-(\varepsilon)$ as a fixed point (a zero).
\emph{Therefore the imaginary part of the displacement map does not vanish
along the open interval $(S_+(\varepsilon),S_-(\varepsilon))$. }
  \item
  The number of the zeros of  the imaginary part of the displacement map,
  along the real analytic curve $\mathcal{H}_\varepsilon^+$ equals the number
  of the zeros of the imaginary part of $d^-_\varepsilon$ along this curve,
  that is to say the number of intersection points of
  $\mathcal{H}_\varepsilon^+$ with $\mathcal{H}_\varepsilon^-$.
\emph{According to  (\ref{hez2}), (\ref{hez4}) and Corollary \ref{simplezero},
this number is one.}
 \end{enumerate}

We conclude that the displacement map can have at most two zeros in the domain
$\mathcal{D}_\varepsilon$, this for all positive $\varepsilon$ smaller than
$\varepsilon_0$ (similar considerations are valid for negative $\varepsilon$).

As we already noted, $d_0=0$ implies $d_1\neq 0$ in the expansion
(\ref{expansion}) and therefore $M_1$ can have at most one simple zero close
to $t=0$. One may wonder, whether two limit cycles can bifurcate from the
two-saddle loop in the case. The somewhat surprising answer is "yes",
as noticed first in \cite{duro06}. The bifurcation of the second  "alien"
limit cycle will be explained in an Appendix. This completes the proof of
Proposition\ref{pr5}. $\Box$

\subsection{The case $M_1=0$}
\label{3.3}
In this section we  suppose that  the Melnikov function $M_1(t)$ vanishes
identically.
The first return map has the form (\ref{melnikov})
where
\begin{equation}
\label{mkt}
M_k(t)=
\int_{\delta(t)}[  \alpha+\beta x+\gamma x^{-1}] ydx,\quad k\geq 2 \quad
\alpha,\beta,\gamma \in \mathbb{R} .
\end{equation}

As we explained, we may suppose that the Bautin ideal is locally principal
at $\lambda_0$ and let $\varepsilon$ be the generator. The deformed vector
field $X_\lambda$ defines a foliation
$$
dH - \sum_{i=1}^\infty \varepsilon^i \omega_i = 0
$$
with first return map
$$
P_{\varepsilon}(h)= h + \varepsilon^k [M_k(h) + O(\varepsilon)],
\quad M_k\neq 0 .
$$
If $\int_{\delta(t) } \omega_1 \not\equiv 0$ then $k=1$ and moreover
$$
M_1(t)= \int_{\delta(t) } \omega_1 .
$$
If, on the other hand, $M_1=0$, then $d \omega_1= c  y dx dy$, where $c$ is
a constant (eventually zero). In general, we shall have
\begin{equation}
\label{chain}
d\omega_1= \dots= d \omega_{d -1} = 0, \;\; d \omega_d= (a+b x+c  y )dx dy
\end{equation}
where
$$
M_d(t)= \int_{\delta(t) } (a+b x+ c y)dx dy .
$$
The case $a^2+b^2\neq 0$ is completely analogous to the case when the first
Melnikov function $M_1$ is not identically zero, and is studied as in
Section \ref{m10}. To the end of the section we consider the case $a=b=0$,
$c\neq0$, in which case the first  non-vanishing Poincar\'e-Pontryagin function
is $M_k$ with suitable $k> d$.

\begin{proposition}
\label{pr6}
If $\gamma\neq 0$, then no limit cycles bifurcate from the two-saddle loop
$\Gamma$.
\end{proposition}

Following the method of the preceding section, we evaluate the number of the
zeros of the displacement map
$$
d^+_\varepsilon - d^-_\varepsilon = (P_\varepsilon - id)\circ d^-_\varepsilon =
 \varepsilon^k M_k(t)+\varepsilon^{k+1}M_{k+1}(t)+ \dots
$$
in the domain $\mathcal{D}_\varepsilon$.

\begin{enumerate}
\item The displacement map,  along the circle $S_R$ is approximated by
$\varepsilon^k M_k(t)$ which has as a leading term $\ln t$ as $\gamma\neq 0$,
see Lemma \ref{asymp}. The increase of the argument of $ \ln t$, and hence of
the displacement map,  along the circle $S_R$ is \emph{close to $0$
 but strictly less than $0$}.
 \item The imaginary part of the displacement map,  along the  interval
 $[S_+(\varepsilon),S_-(\varepsilon)]$ equals  the imaginary part of
$d^-_\varepsilon(t)$. Its zeros equal the number of intersection points
of $\mathcal{H}_\varepsilon^-$ with the real axes, which amounts to solve
$h^\varepsilon_{\delta_-}(z)=z$, see
(\ref{hez}). Zeros of $h^\varepsilon_{\delta_-}-id$ correspond to complex
limit cycles (except the origin $S_-$). Their number is the \emph{cyclicity
of the saddle point}. We have
$$
h^\varepsilon_{\delta_-}(z)=z + \varepsilon^d M^-_d(t)+ \dots,
\quad a, b, c \in \mathbb{R}
$$
where
$$
M_d^-(t) = \int_{\delta_-(t)} \omega_d, \;\; d \omega_d = c y dx dy , c\neq 0 .
$$
 Lemma \ref{asymp} implies $ \int_{\delta_-(t)} y^2 dx  =  \pm 2 \pi i t$, and
 hence the cyclicity of the saddle point is zero.
We conclude that the imaginary part of the displacement map does not vanish
along the  interval $[S_+(\varepsilon),S_-(\varepsilon))$.
  \item The number of the zeros of  the imaginary part of the displacement map,
  along the real analytic curve $\mathcal{H}_\varepsilon^+$ equals the number of
  zeros of the imaginary part of $d^-_\varepsilon$ along this curve, that is to
  say the number of intersection points of $\mathcal{H}_\varepsilon^+$ with
$\mathcal{H}_\varepsilon^-$. According to  (\ref{hez2}) we need the expansion
of $h^\varepsilon_{\delta_\pm}(z)-z$. The monodromy of the first return map
$P_\varepsilon(e^{2\pi i }t)- P_\varepsilon(t)$, equals the holonomy
$h_{\delta_-}^\varepsilon\circ h_{\delta_+}^\varepsilon(z)$, where
$z$ is a different chart close to $t$, $z=t+O(\varepsilon)$. Therefore, if

$$
P_\varepsilon(t)=t + \varepsilon^k (\ln t \int_{\delta_+(t)+ \delta_-(t)}
[\alpha+\beta x+\gamma x^{-1}] ydx + h.f.) + O(\varepsilon^{k+1})
$$
then
$$
h^\varepsilon_{\delta_-}\circ h^\varepsilon_{\delta_+}(z)= 2 \pi i
\varepsilon^k \int_{\delta_+(t)+ \delta_-(t)}[  \alpha+\beta x+
\gamma x^{-1}] ydx + O(\varepsilon^{k+1}) .
$$
The notation $O(\varepsilon^{k+1})$ has as usual an appropriate meaning.
It represents a function which, for a fixed $z$ or $t$, is bounded by a
function of the type $ O(|\varepsilon|^{k+1})$. Finally, "h.f." stays for
a function, holomorphic in $t$. As the leading term of $P_\varepsilon(t)$
is $\ln t$ multiplied by a non-zero constant, then the above formula shows
that the leading term of the holonomy map is a non-zero constant
$$h^\varepsilon_{\delta_-}\circ h^\varepsilon_{\delta_+}(z) = \varepsilon^k
(c+ \dots) +  O(\varepsilon^{k+1}),\; c\neq 0 .$$
The conclusion is that the imaginary part of the displacement map has no
zeros along the real analytic curve $\mathcal{H}_\varepsilon^+$.
 \end{enumerate}
Summing up the above information, we conclude that  the displacement map has
no zeros in the domain $\mathcal{D}_\varepsilon$. Proposition \ref{pr6} is
proved. $\Box$
\begin{proposition}
\label{pr7}
If $\gamma =  0$, but $\alpha\neq 0$, then at most two limit cycles bifurcate
from the two-saddle loop $\Gamma$.
\end{proposition}
{\bf Proof.}
The condition $\alpha\neq 0$ is equivalent to the condition $d_1\neq 0$ in the
expansion of the first non-vanishing Melnikov function
$$
M_k(t)= d_0+ d_1 t \ln t + d_2 t + d_3 t^2 \ln t + \dots
$$
 \begin{enumerate}
 \item
 The displacement map,  along the circle $S_R$ is approximated by
 $\varepsilon^k M_k(t)$ which has as a leading term
 either a constant, to $t \ln t$. In both cases
 the increase of the argument  of the displacement map,  along the circle
 $S_R$ is  strictly less than $2\pi$.
 \item The imaginary part of the displacement map,  along the  interval
 $[S_+(\varepsilon),S_-(\varepsilon)]$ equals  the imaginary part of
$d^-_\varepsilon(t)$. As in the preceding proposition, we get that the
imaginary part of the displacement map does not vanish along the interval
$[S_+(\varepsilon),S_-(\varepsilon))$.
\item The number of the zeros of  the imaginary part of the displacement map,
along the real analytic curve $\mathcal{H}_\varepsilon^+$, equals the number
of intersection points of this curve with $\mathcal{H}_\varepsilon^-$. It is
bounded by the cyclicity of
$$
d_1 t + d_3 t^2 + \dots
$$
that is to say by one. This implies the statement of Proposition \ref{pr7}.
 $\Box$

 \end{enumerate}
\begin{proposition}
\label{pr8}
If $\gamma= \alpha =0$, but $\beta\neq 0$,
then at most three limit cycles bifurcate from the two-saddle loop $\Gamma$.
\end{proposition}
The condition $\alpha= 0$ but $\beta\neq 0$ implies $d_1= 0$, $d_3\neq 0$,
$d_0\neq 0$ in the expansion of the first non-vanishing Melnikov function
$$
M_k(t)= d_0+ d_1 t \ln t + d_2 t + d_3 t^2 \ln t + \dots
$$
Repeating the preceding arguments, we obtain a bound of three limit cycles
(possibly complex). $\Box$

\section{Global results}
\label{globalsection}
Let $H(x,y)$ be a real cubic polynomial, such that $X_H$ has a non-degenerate
two-saddle loop $\Gamma$ as on Figure \ref{fig1}. Denote by $\Pi$ the period
annulus surrounded by $\Gamma$, and by $\bar{\Pi}= \Pi\cup \Gamma$ its closure.
Theorem \ref{main1} can be generalized as follows
\begin{theorem}
\label{semiglobal}
The cyclicity of the closed period annulus $\bar{\Pi}$ under an arbitrary
quadratic deformation, is less then or equal to three.
\end{theorem}
Let $X_\varepsilon$ be a one-parameter family of plane quadratic vector
fields, depending analytically on a real parameter $\varepsilon$, and such
that $X_0=X_H$ is a Hamiltonian vector field having a non-degenerate two-saddle
loop $\Gamma$ as above.

\begin{theorem}
\label{semilocal}
If the first Melnikov function is not identically zero, and
\begin{itemize}
\item $M_1(0)\neq 0$, then no limit cycles bifurcate from $\Gamma$ and at most
one limit cycle bifurcates from the closed period annulus $\bar{\Pi}$;
\item
$M_1(0)=0$, then at most two limit cycles bifurcate from the two-saddle loop
$\Gamma$ and no limit cycles bifurcate from the open period annulus $\Pi$.
\end{itemize}
If the first non-vanishing Melnikov function $M_k$, $k\geq 2$ is as in
(\ref{mkt}), and
\begin{itemize}
\item
$\gamma \neq 0$, then no limit cycles bifurcate from $\Gamma$ and at most two
limit cycles bifurcate from the closed period annulus $\bar{\Pi}$;
\item $\gamma = 0$ and $M_k(0)=0$, then at most two limit cycles bifurcate from
the two-saddle loop $\Gamma$ and no limit cycles bifurcate from the open period
annulus $\Pi$;
\item  $\gamma = 0$, $\alpha\neq 0$ and $M_k(0)\neq 0$, then no limit cycles
bifurcate from $\Gamma$ and at most one limit cycle bifurcates from the
closed period annulus $\bar{\Pi}$;
\item $\gamma=\alpha=0$ and $\beta\neq 0$, then no limit cycles bifurcate from
the open period annulus $\Pi$, and at most three limit cycles bifurcate from
the two-saddle loop $\Gamma$.
\end{itemize}
\end{theorem}

Let $H(x,y)$ be a real cubic polynomial with four distinct (real or complex)
critical points, but only three distinct critical values.
Let $X_H$ be the corresponding quadratic Hamiltonian vector field
(\ref{hamiltonian}).
\begin{theorem}
\label{global}
There is a neighborhood
${\cal U}$
of $X_H$
in the space of all quadratic
vector fields, such that any $X
\in
{\cal U}$
has at most three limit cycles.
\end{theorem}

Theorem \ref{global} is the analogue of  \cite[Theorem 1]{gavr01}, \cite[Theorem 2]{hoil94b} , where it
is shown that for a cubic Hamiltonian $H(x,y)$ with four distinct critical
values, the exact upper bound for the number of the limit cycles of any
sufficiently close quadratic system, is two. Let us explain in brief which
$X_H$ Theorem \ref{global} concerns. By using the normal form for cubic
Hamiltonians with a center from \cite{hoil94b},
$$H(x,y)=\frac{x^2+y^2}{2}-\frac{x^3}{3}+axy^2+\frac{b}{3}y^3, \;\;
-\frac12\leq a\leq 1,\;\; 0\leq b\leq (1-a)\sqrt{1+2a},$$
one can easily verify that the level value corresponding to a critical point
$(x_0,y_0)$ is $H(x_0,y_0)=\frac16(x_0^2+y_0^2)$. Then, for the generic
Hamiltonians (corresponding to internal points $(a,b)$ of the domain of
parameters) there are either four distinct critical levels or three distinct
critical points in the finite plane and Theorem \ref{global} does not
concern them. For the degenerate Hamiltonians (corresponding to points
from the boundary of the domain of parameters), there are four distinct
critical points with three distinct critical values if and only if
$(a,b)\neq (-\frac12,0), (-\frac13,0), (0,0), (1,0),
(\frac12,\sqrt{\frac12})$. Therefore, in the normal form (\ref{ham}),
Theorem \ref{global} concerns all $a\in {\mathbb R}$ except $a=-1, 0,2,3$.

{\bf Conjecture.} {\it The exact upper bound for the number of limit cycles
in Theorem \ref{main1}, Theorem \ref{semiglobal} and Theorem \ref{global}
is two.}\\ \\

\begin{figure}
\begin{center}
\def\svgwidth{10cm}
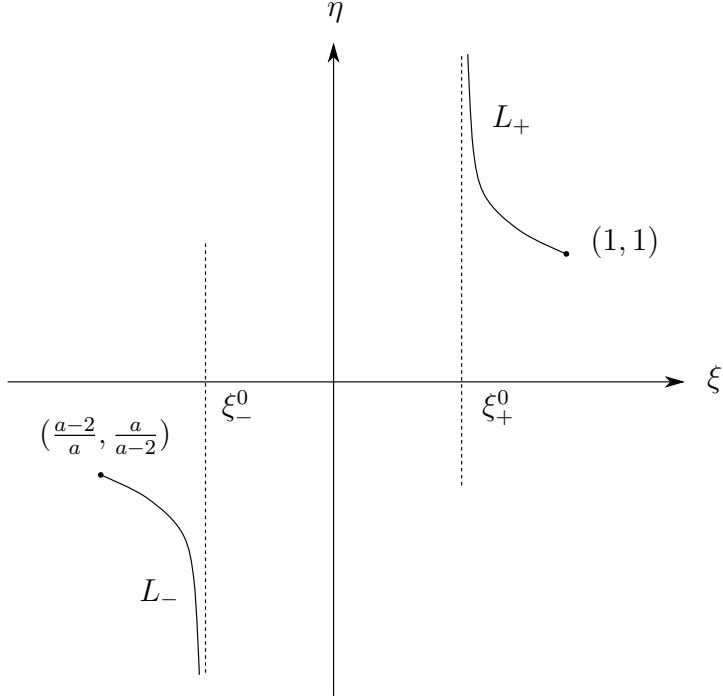

\end{center}
\label{centroid}
\caption{The curves $L_+$ and $L_-$}
\end{figure}
{\bf Proof of Theorems \ref{semiglobal},\ref{semilocal},\ref{global}.}
For the saddle-loop cases (that is $a\not\in[-1,2]$) in Theorem \ref{global},
it is well known that at most two limit cycles can bifurcate from the closed
period annulus \cite{gail00, cly02}. Below we are going to apply the results
just established to handle the two-saddle loop cases $a\in (-1,2)$.
The proofs will follow from a careful comparison of the statements in the
preceding section and the available results on the cyclicity of open period
annuli of quadratic Hamiltonian systems, see \cite{zhzh01,ilie96,cly02}.

Using the notations of Section \ref{abelianintegrals},
denote by $\Sigma_+=[a-3,0)$ the semi-open interval with respect to $t$
corresponding to the period annulus surrounding the center $C_+$ at $(1,0)$.
When there is a second center $C_-$ at $(\frac{a-2}{a},0)$
which happens for $0<a<2$, we shall denote the related interval by
$\Sigma_-=(0,\frac{(a+1)(a-2)^2}{a^2}]$. Consider the respective Melnikov
function(s)
$$M_k(t)=\alpha J_0(t)+\beta J_1(t)+\gamma J_{-1}(t),\qquad t\in \Sigma_\pm.$$
Next, define the planar curve(s)
$$L_\pm=\left\{(\xi_\pm(t),\eta_\pm(t))=\left(\frac{J_1(t)}{J_0(t)},
\frac{J_{-1}(t)}{J_0(t)}\right):\quad t\in \Sigma_\pm\right\}.$$
The properties of the curves $L_\pm$ are well known, see
\cite{zhzh01}, \cite{ilie96} and \cite{cly02} for the hyperbolic,
the parabolic and the elliptic cases. Namely (see Figure \ref{centroid}),

\vspace{2ex}\noindent
1) $\xi_+(t)$ is decreasing, $\eta_+(t)$ is increasing and $L_+$ is a convex
   curve. $L_+$ begins at point $(1,1)$ and has a vertical asymptote
   $\xi=\xi_+(-0)=c_0/b_0$ as $t\to -0$.

\vspace{2ex}\noindent
2) If $L_-$ exists, then $\xi_-(t)$ is decreasing, $\eta_-(t)$ is increasing
and $L_-$ is a concave curve. $L_-$ ends at point $(\frac{a-2}{a},\frac{a}{a-2})$
and has a vertical asymptote $\xi=\xi_-(+0)$ as $t\to +0$.

\vspace{2ex}\noindent
3) The number of limit cycles born from periodic orbits equals the number of
the intersections (counted with multiplicities) between the straight line
$\alpha+\beta\xi+\gamma\eta=0$ and the curve $L_+$ (both curves $L_\pm$
in the elliptic case).

\vspace{2ex}\noindent
4) If $P_*$ is intersection point corresponding to $t=t_*$, then the
related limit cycle approaches the oval $H(x,y)=t_*$ as
$\varepsilon\to 0$.

\vspace{2ex}\noindent
Now, if $\gamma\neq 0$, then by Proposition \ref{pr6} above, there are no limit
cycles produced by the double loop(s). On the other hand, any line has at most
two intersection points with $L_\pm$. Two is the total upper bound of the
number of limit cycles produced under the perturbation.

Next, if $\gamma=0$, then by Proposition \ref{nocycles}, a necessary condition
for the
bifurcation of limit cycles from the double loop(s) is
$\alpha J_0(0)+\beta J_1(0)=0$. It is easy to see that limit cycles cannot
bifurcate simultaneously from both two-saddle loops existing when $a\in (0,2)$.
Indeed, the system
$$\alpha J_0(-0)+\beta J_1(-0)=\alpha J_0(+0)+\beta J_1(+0)=0$$
implies $\alpha=\beta=0$. This is because the system is equivalent to
$$\alpha+\xi_+(-0)\beta=\alpha +\xi_-(+0)\beta=0\;\;\mbox{\rm and}\;\;
\xi_-(+0)<0<\xi_+(-0).$$
Therefore, if $\gamma=0$ but $\alpha\neq 0$, then by Proposition \ref{pr7}
above, there are
at most two limit cycles produced by the double loop(s). On the other hand,
any line $\alpha+\beta\xi=0$ has at most one intersection point with $L_\pm$.
Moreover, if such a point exists, no limit cycles are produced by the double
loop(s), according to  Proposition \ref{nocycles}. Again, two is the total
upper bound of the number of limit cycles produced under the perturbation.

If $\gamma=\alpha=0$ but $\beta\neq 0$, then by Proposition \ref{pr8} above,
there are at most three limit cycles produced by the double loop(s). On the
other hand, the line $\xi=0$ has no intersection points with $L_\pm$.
Hence, three is the total upper bound of the number of limit cycles produced
under the perturbation. $\Box$

\section{Appendix : alien limit cycles in quadratic systems}
\begin{figure}
\begin{center}
\def\svgwidth{10cm}
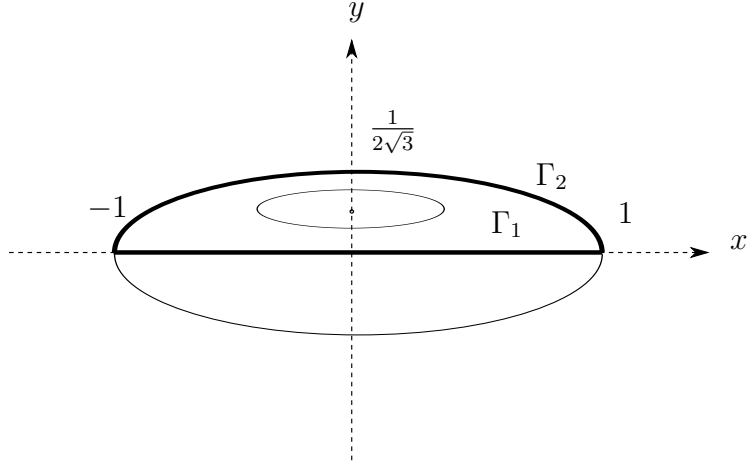
\end{center}
\label{app}
\caption{The two-saddle loop $\Gamma_u=\Gamma_1\cup \Gamma_2 \subset \{H=0\}.$ }
\end{figure}

Consider, using the notations of \cite{ldcr09}, the perturbed quadratic
Hamiltonian system
\begin{equation}
\label{perturbed1}
X_{\mu,\varepsilon}: \left\{
\begin{array}{lcl}
   \dot{x}   &= &H_y \\
       \dot{y}   &= & -H_x - \varepsilon  P
\end{array} \right.
\end{equation}
where
\begin{equation}
\label{ham2}
H= y (x^2+ \frac {1}{12} y^2 -1), \quad P(x,y,\mu) = (16+c x- \pi \sqrt{3} y)y
+ \mu_1 + \mu_2 y , \;\;
\end{equation}
$\varepsilon, \mu_1, \mu_2$ are sufficiently small real numbers, and $c$ is a
real constant bigger than 16.
Denote the upper two-saddle loop of the non perturbed system ($\varepsilon = 0$)
by $\Gamma_u=\Gamma_1\cup \Gamma_2$, where $\Gamma_1$ is the segment
$\{ (x,y): -1\leq x \leq 1, y=0 \}$ and $\Gamma_2$ is the half-ellipse $\{ (x,y)
: x^2+ \frac {1}{12} y^2 =1, y \geq 0 \}$, see Fig. \ref{app}.
Let
$$
\{\gamma(h)\}_h \subset \{(x,y) \in \R^2 : H(x,y)=h\}
$$
be the continuous family of ovals, contained in the two-saddle loop $\Gamma_u$,
parameterized by $h\in (-4/3,0)$.
The first return map  of $X_{\mu,\varepsilon}$ takes the form
$$
h \mapsto h + \varepsilon \int_{\gamma(h)} P(x,y,\mu) dx + O(\varepsilon^2)
$$
where $\int_{\gamma(h)} P(x,y,\mu) dx$  is the first Poincar\'e-Pontryagin function associated to $X_{\mu,\varepsilon}$.
We have
$$
\int_{\gamma(h)} P(x,y,\mu) dx = d_0(\mu)+ d_1(\mu) h \log(h) + O(h)
$$
see   (\ref{expansion}). It is straightforward to check that 
$d(0)=0$ and by Proposition \ref{d0d1} then we get $d_1(0)\neq 0$. It follows that for sufficiently small $\|\mu\|$, $|h|$, $h<0$, the Poincar\'e-Pontryagin function $\int_{\gamma(h)} P(x,y,\mu) dx$ has at most one zero. The purpose of this Appendix is to show that
the number of the limit cycles, which bifurcate from $\Gamma_u$, exceeds the number of the zeros of 
$\int_{\gamma(h)} P(x,y,\mu) dx$
near $h=0$. The "missing" second limit cycle, which does not correspond to a zero is an
"alien" limit cycle.This is  a new unexpected phenomenon in the bifurcation
theory of vector fields, discovered recently by Caubergh, Dumortier and
Roussarie \cite{cdr05, duro06}. In contrast to the preceding examples
 \cite{cdr07, ldcr09, cdl10,codupr13}) the system which we consider is quadratic.
\begin{proposition}
\label{cycl2}
 \textit{The cyclicity $Cycl(\Gamma_u,X_{\mu,\varepsilon})$ of the two-loop
 $\Gamma_u$ with respect to the deformed vector field $X_{\mu,\varepsilon}$
 is two}.
 \end{proposition}
Note that, according to Proposition \ref{pr5}, the cyclicity
$Cycl(\Gamma_u,X_{\mu,\varepsilon})$ is at most two.\\
{\bf Proof of Proposition  \ref{cycl2}.}
We shall follow closely \cite[section 6.2.]{duro06}. The traces  $\sigma_{1,2}$
of the vector field $X_{\mu,\varepsilon}$ at the saddle points determine its
"stability".
%
%
\begin{figure}
\begin{center}
\resizebox{15cm}{!}
{ \includegraphics{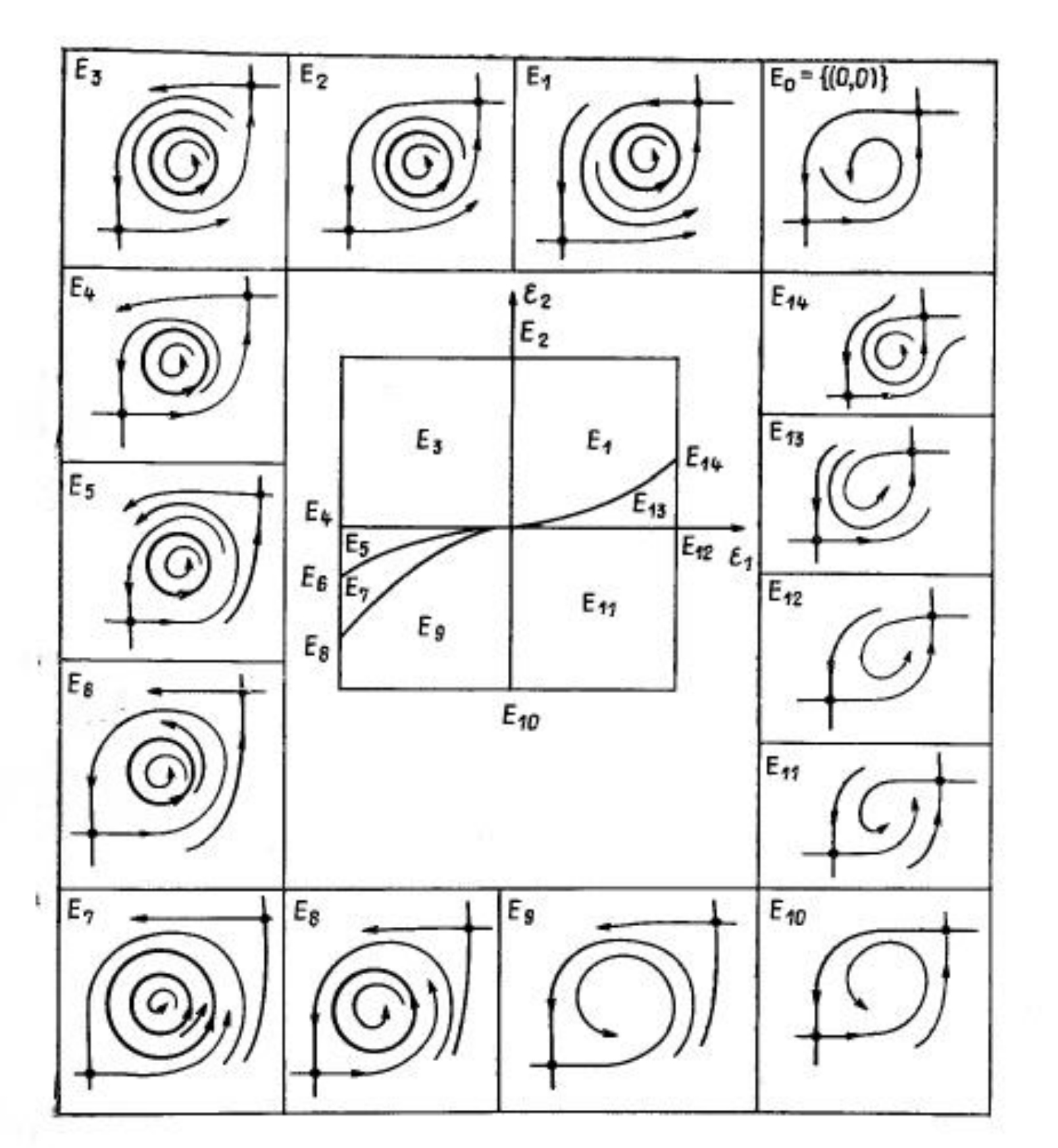} }
\end{center}
\caption{Bifurcation diagram of generic two-parameter deformations of vector
fields, containing a two-saddle loop. In the domain $E_7$ the system has two
limit cycles.}
\label{figr}
\end{figure}
As the coordinates of the saddle points satisfy
$$
x= \pm 1 + O(\varepsilon), \;\; y= O(\varepsilon)
$$
then
for the traces $\sigma_{1,2}$ at the saddle points $s_1, s_2$ we get
$$
\sigma_1(\varepsilon,\mu)= (-16+c- \mu_2 ) \varepsilon + O(\varepsilon^2)
$$
$$
\sigma_2(\varepsilon,\mu)= (-16-c- \mu_2 ) \varepsilon + O(\varepsilon^2)
$$

For small $\varepsilon$ and a general perturbation, the connections
$\Gamma_{1,2}$ will be broken. The distance between the two branches (stable
and unstable separatrix) of the broken connection can be measured on a segment,
transverse to $\Gamma_1$ or $\Gamma_2$. Let us denote these distances
(or shift functions) by $b_{1,2}$. It is well known that the shift functions
are analytic functions in $\varepsilon, \mu$, and if we use the restriction of
$H$ to the transverse segments as a local parameter $h$, then
\begin{equation}
\label{shift}
b_{i}(\varepsilon,\mu) = \varepsilon \int_{\Gamma_i} \omega_\mu +
O(\varepsilon^2), \quad i=1,2.
\end{equation}
With the notations above we compute
$$
\int_{\Gamma_2} y dx = - \pi \sqrt{3}, \quad \int_{\Gamma_2} y^2 dx = -16
$$
and therefore
$$
\int_{\Gamma_2} \omega_\mu = -2 \mu_1 - \pi \sqrt{3} \mu_2, \qquad
\int_{\Gamma_1} \omega_\mu = 2 \mu_1 .
$$
It is immediately seen that
\begin{itemize}
\item for every sufficiently small $\varepsilon\neq 0$ and $\| \mu\|$, the
traces $\sigma_1, \sigma_2$ are non-zero and have opposite signs;
\item for every sufficiently small $\varepsilon\neq 0$ and $\| \mu\|$
$$
\det \left(\begin{array}{cc}\frac{\partial b_1}{\partial \mu_1} &
\frac{\partial b_1}{\partial \mu_2}  \\\frac{\partial b_2}{\partial \mu_1}  &
\frac{\partial b_2}{\partial \mu_2} \end{array}\right) \neq 0.
$$
\end{itemize}
Under these conditions, the bifurcation diagram of limit cycles near the
double loop $\Gamma_1 \cup \Gamma_2$ was computed by  Dumortier, Roussarie
and  Sotomayor \cite{drs91}, see
\cite[fig. 5]{duro06}. It follows that the cyclicity of the two loop $\Gamma$
under the quadratic perturbation (\ref{perturbed1}) is two.$\Box$\\
{\bf Remark.}
An alternative proof of  Proposition \ref{cycl2} can also be obtained from the
classical Roitenberg Theorem, see  \cite[Theorem 2, fig. 40a]{aais94}, which
is illustrated on Fig. \ref{figr}. Namely, as the deformation
(\ref{perturbed1}) depends on three parameters, then there is a one-parameter
induced deformation
\begin{equation}
\label{mue}
\mu_1=\mu_1(\varepsilon) = O(\varepsilon),\quad \mu_2=\mu_2(\varepsilon)
= O(\varepsilon)
\end{equation}
such that the two connections $\Gamma_1$ and $\Gamma_2$ persist for all
sufficiently small $\varepsilon$. This one-parameter deformation is not in an
integrable direction at a first order in $\varepsilon$, in the sense that the
corresponding first Melnikov function $M_1(h,\mu)|_{\mu=0}$ is not identically
zero. One easily verifies that this implies the genericity assumptions of
 \cite[Theorem 2]{aais94}. Thus, making an additional deformation in a
 direction transversal to the curve (\ref{mue}), we get the bifurcation diagram
 of Roitenberg shown on Fig. \ref{figr}. This diagram  is a two dimensional
 section $\{ \varepsilon= const\}$ of the three-dimensional diagram
 \cite[fig.5]{duro06}.


\vspace{2ex} \noindent {\bf Acknowledgments.}
We are obliged to V. Roitenberg who gave us a permission to reproduce
Fig. \ref{figr}.
Part of this work has been done while the second author visited the
University of Toulouse. He is very grateful for kind hospitality.

\bibliography{bibliography}
\end{document}